\documentclass
{article}
\usepackage[T2A]{fontenc}
\usepackage[cp1251]{inputenc}
\usepackage[english,russian]{babel}
\usepackage[tbtags]{amsmath}
\usepackage{amsfonts,amssymb,mathrsfs,amscd}
\usepackage{wasysym}
\usepackage{verbatim}
\usepackage{eucal}
\usepackage{upgreek}
\usepackage{graphicx}
\usepackage{enumerate}
\let\No=\textnumero

\usepackage[hyper]{izv2e}
\JournalName{
}

\numberwithin{equation}{section}

\theoremstyle{plain}
\newtheorem{theorem}{Теорема}

\newtheorem{lemma}{Лемма}
\newtheorem{propos}{Предложение}

\newtheorem{thN}{Теорема Р. Неванлинны}

\newtheorem{lemGS}{Лемма Гришина\,--\,Содина о малых интервалах}

\theoremstyle{definition}
\newtheorem{definition}{Определение}
\newtheorem{proof}{Доказательство}
\newtheorem{remark}{Замечание}
\newtheorem{example}{Пример}

\renewcommand{\leq}{\leqslant} 
\renewcommand{\geq}{\geqslant}
\newcommand{\RR}{\mathbb{R}} 
\newcommand{\CC}{\mathbb{C}} 
\newcommand{\NN}{\mathbb{N}}

\newcommand{\rad}{\text{\tiny\rm rd}}
\renewcommand{\Re}{{\rm Re\,}}

\DeclareMathOperator{\supp}{{\sf supp}}
\DeclareMathOperator{\dd}{\,{\mathrm  d\!}}

\DeclareMathOperator{\mes}{mes}

\begin{document} 
\title{Мероморфные функции и разности субгармонических функций в интегралах и разностная характеристика Неванлинны. I. Радиальные максимальные характеристики роста}

\author[B.\,N.~Khabibullin]{Б.\,Н.~Хабибуллин}
\address{Башкирский государственный университет}
\email{khabib-bulat@mail.ru}

\date{}
\udk{517.547.26 : 517.547.28 : 517.574}

 \maketitle

\begin{fulltext}

\begin{abstract} Пусть  $f$ --- мероморфная функция на комплексной плоскости $\mathbb C$ с функцией максимума её модуля    $M(r,f)$ на окружностях с центром в нуле радиуса $r$. Ряд классических, известных и широко используемых результатов позволяют оценить сверху интегралы от положительной части логарифма $\ln^+M(t,f)$ по подмножествам $E$ на отрезках $[0,r]$ через характеристику Неванлинны $T(r,f)$ и  линейную лебегову меру множества $E$. Статья даёт  аналогичные оценки для интегралов Лебега\,--\,Стилтьеса от $\ln^+M(t,f)$ по возрастающей функции интегрирования $m$.   Основная часть изложения ведётся сразу для разностей субгармонических  функций в кругах с центром в нуле, или $\delta$-субгармонических функций. Единственное условие  в основной теореме  ---  модуль непрерывности функции интегрирования $m$ удовлетворяет условию Дини. Это условие в некотором смысле и необходимо. Таким образом, первая часть работы в определённой степени завершает в общей форме исследования по верхним оценкам интегралов от радиальных максимальных характеристик роста  произвольных мероморфных и $\delta$-субгармонических функций  через характеристику Неванлинны с её версиями и через характеристики функции интегрирования $m$. 

Библиография:  17 названий 

Ключевые слова: мероморфная функция, $\delta$-субгармоническая функция, характеристика Неванлинны, мера Рисса, модуль непрерывности 

\end{abstract}

\markright{Мероморфные функции и разности субгармонических функций  \dots}


\section{Введение}

Как обычно, $\NN:=\{1,2,\dots\}$ ---   множество 
{\it натуральных чисел,\/}   $\CC$ --- {\it комплексная плоскость\/} с  {\it вещественной осью\/}   $\RR$ и {\it положительной полуосью\/}  $\RR^+:=\{x\in \RR \,|\, x\geq 0\}$. Одноточечные множества записываем без фигурных скобок, если это не вызывает разночтений. Так, $\overline \RR:=\RR\cup \pm\infty$ и $\overline \RR^+:=\RR^+\cup+\infty$   --- {\it расширенные\/} соответственно  вещественная ось  и положительная полуось с обычным {\it модулем\/} $|\cdot|$, как и для $\CC$,  и с  $|\pm \infty|:=+\infty$, а также $\sup \varnothing :=-\infty=\inf \RR$, $\inf \varnothing:=+\infty=\sup \RR $ для {\it пустого множества\/} $\varnothing$. 
Для  $+\infty\in \overline \RR$ используем и обозначение $\infty$ без знака $+$.
{\it Интервал $I$ на\/ $\overline \RR$} --- это {\it связное\/} подмножество в $\overline \RR$ с {\it левым концом\/}
$\inf I\in \overline \RR$ и с {\it правым концом\/} $\sup I\in \overline \RR$. 
{\it Отрезок\/} с {\it концами\/} $a\leq b\in \overline \RR$ --- это интервал $[a,b]:=\bigl\{x\in \overline \RR\bigm| a\leq x\leq b\bigr\}$. Кроме того,  $[a,b):=[a,b]\setminus b$ (соответственно  $(a,b]:=[a,b]\setminus a$) --- {\it открытый справа}
(соответственно {\it слева}) и {\it замкнутый слева\/} (соответственно {\it справа}) интервал, 
$(a,b):=[a,b)\cap (a,b]$ --- {\it открытый интервал.\/} 

Через $x^+:=\sup \{0,x\}\in \overline \RR^+$ обозначаем {\it положительную часть\/} элемента  $x\in \overline \RR$, а $x^-:=(-x)^+\in \overline \RR^+$ --- его {\it отрицательная часть\/}.    Вообще всюду далее
{\it положительность\/} --- это $\geq 0$, а {\it отрицательность\/} --- это $\leq 0$. Если $x\in \overline \RR^+\setminus 0$, то $x$ --- {\it строго\/} положительный элемент, а если $x\in -\overline\RR^+\setminus 0$, то {\it строго\/} отрицательный элемент.  

Для {\it расширенной числовой функции\/} $f\colon X\to \overline \RR$ через 
$f^+\colon x\underset{\text{\tiny $x\in X$}}{\longmapsto} \bigl(f(x)\bigr)^+$
обозначаем её {\it положительную часть\/}, а $f^{-}:=(-f)^+$ --- её {\it отрицательную  часть.\/}
При $X\subset \RR$ эта функция $f$ {\it возрастающая,\/} если  для любых  $x_1, x_2\in X$ из $x_1\leq x_2$ следует $f(x_1)\leq f(x_2)$, а возрастающая функция $f$ {\it строго\/} возрастающая, если она инъективна.   
Функция $f$ {\it (строго) убывающая,\/} если противоположная функция $-f$ (строго) возрастающая.
{\it Числовая функция\/} $f\colon I\to \RR$ на интервале $I\subset \overline \RR^+$ {\it выпукла относительно\/} логарифма $\ln$, если суперпозиция $f\circ \exp$ выпукла на интервале $\ln I:=\bigl\{x\in \overline \RR\bigm| e^x\in I\bigr\}$.

\subsection{Предшествующие  результаты}\label{s10}

Как исходная точка основных результатов настоящей статьи может рассматриваться одна  из классических теорем Рольфа Неванлинны  \cite[pp. 24--27]{RNevanlinna}, на которую имеется ссылка   
в классической монографии А.\,А.~Гольдберга и И.\,В.~Островского  \cite[Комментарии, К главе I]{GO}, 
где она также и доказана \cite[гл.~I, теорема 7.2]{GO}. Приведём её здесь в обозначениях и определениях, а также в формулировке  из \cite{GO}.

Пусть $f$ --- мероморфная функция на {\it комплексной плоскости\/} $\CC$ с {\it функцией максимума модуля\/}    
 на $\RR^+$, определённой  как
\begin{subequations}\label{TN}
	\begin{align}
	M(r,f)&:=\sup\Bigl\{\bigl|f(z)\bigr| \Bigm| |z|=r\Bigr\} \quad\text{для  $r\in \RR^+$}
	\tag{\ref{TN}M}\label{{TN}M}\\
	\intertext{и с {\it характеристикой Неванлинны}}
	T(r, f)&\underset{r\in\RR^+}{:=}m(r,f)+N(r,f),  \quad\text{где}
	\tag{\ref{TN}T}\label{{TN}T}
	\\
	m(r,f)&\underset{r\in\RR^+}{:=}\frac{1}{2\pi}\int_0^{2\pi} \ln^+\bigl|f(re^{i\varphi})\bigr| \dd \varphi , 
\tag{\ref{TN}m}\label{{TN}m}\\
N(r,f)&\underset{r\in\RR^+}{:=}\int_{0}^{r}\frac{n(t,f)-n(0,f)}{t}\dd t+n(0,f)\ln r, 
	\tag{\ref{TN}N}\label{{TN}N}
	\end{align}
\end{subequations} 
а   $n(r,f)$ --- число полюсов функции  $f$ в  замкнутом круге радиуса $r\in \RR^+$ с центром в нуле, подсчитанное с учётом кратности полюсов. Элементарные их свойства по $r$  --- это {\it положительность\/} $M$, $n$ и $m$, 
{\it возрастание\/}  $n$, $N$ и  $T$, {\it выпуклость\/} $N$ и $T$ {\it относительно\/} $\ln$ и {\it непрерывность\/}  $N$, $T$ и $m$ при $r>0$. 

\begin{thN}[{\rm  \cite[гл. I, теорема 7.2]{GO}}] 
Пусть   $f(z)$  --- мероморфная функция, $k$ --- некоторое число, $k>1$. Тогда справедливо неравенство 
\begin{equation}\label{eN}
	\frac{1}{r}\int_0^r\ln^+ M(t,f)\dd t\leq c(k)T(kr,f), 
	\end{equation}
где постоянная  $c(k)>1$ зависит только от  $k$.   
\end{thN}
Обсуждение в \cite[Введение, 1.1]{Kha20} показывает, что  в \eqref{eN} невозможно обойтись
без ограничения вида $r\geq r_0>0$ для  некоторого фиксированного $r_0>0$, а постоянная $c(k)$ всё же зависит от выбора числа $r_0>0$,  если не добавлять в правую часть 
\eqref{eN} некоторое слагаемое логарифмического порядка роста. 

Оценки интегралов по малым подмножествам на окружности или на луче также широко использовались в теории целых и мероморфных функций.  Исходная  точка таких исследований  --- 
лемма А.~Эдрея и В.~Фукса о малых дугах 
\cite[{\bf 2}, лемма III, {\bf 9}]{EF}, \cite[теорема 7.3]{GO}, которая нашла важные  применения, отражённые, например, в  
\cite[гл.~I, теорем  7.4]{GO}.  Вариант леммы Эдрея\,--\,Фукса о малых дугах --- лемма А.\,Ф.~Гришина и М.\,Л.~Содина  о малых интервалах. Она представлена в \cite[лемма 3.1]{GrS} с применениями \cite[лемма 3.2, теорема 3.1]{GrS},  а  её доказательство, как отмечают авторы, дословно повторяет доказательство леммы  Эдрея\,--\,Фукса о малых дугах. 
Следуя  \cite{GO}, \cite{GrS}, используем распространённое обозначение $\mes$ для {\it линейной меры Лебега\/} на  $\RR$. Подмножество в $\RR$ или расширенную числовую функцию на подмножестве из $\RR$ называем измеримыми, если они $\mes$-измеримы,  а  $\mes E:=\mes (E)$ для измеримого $E\subset \RR$.     

\begin{lemGS}[{\rm \cite[лемма 3.1]{GrS}}] 
Пусть  $f$ --- мероморфная функция, $E\subset [1,+\infty)$ измеримо и  $E(r):=E\cap [1,r)$. Тогда
	\begin{equation}\label{eGS}
	\frac{1}{r}\int_{E(r)}\ln^+M(t,f)\dd t\leq 	C\frac{k}{k-1}\Bigl(\frac{\mes E(r)}{r}\ln 	
\frac{2r}{\mes E(r)} \Bigr)T(kr,f),
	\end{equation}
где  $C$ --- абсолютная постоянная.  
\end{lemGS}

Версия леммы Гришина\,--\,Содина {\it для субгармонических функций конечного порядка\/}  доказана в совместной работе А.\,Ф.~Гришина и Т.\,И.~Малютиной  \cite{GriMal05}, посвящённой ряду вопросов теории роста субгармонических функций конечного порядка \cite{Rans}, \cite{HK}. Теорема Гришина\,--\,Малютиной о малых интервалах недавно  была распространена на произвольные субгармонические функции на $\CC$ в совместной статье Л.\,А.~Габдрахмановой и автора   \cite{GabKha20},  
из результатов которой она легко следует \cite[вывод теоремы Гришина\,--\,Малютиной]{GabKha20}.  Ещё один подобный  более общий результат   установлен в  \cite[основная теорема]{Kha20}. Его формулировка и новый вывод будут даны во второй части работы. Отметим лишь, что  он  содержит в себе как частный случай все приведённые и обсуждённые выше результаты, за исключением леммы Эдрея\,--\,Фукса о малых дугах.

\subsection{Некоторые определения и обозначения} 
Через 
$$D(r):=\bigl\{z \in \CC \bigm| |z|< r\bigr\}, \quad
\overline  D(r):=\bigl\{z \in \CC \bigm| |z|\leq r\bigr\}, \quad 
\partial \overline  D(r):=\overline  D(r)
\setminus  D(r)
$$ 
обозначаем соответственно {\it открытый\/} и  {\it замкнутый круги,\/} а также  {\it  окружность\/}  в $\CC$ {\it радиуса $r\in \overline \RR^+$ с центром в нуле;\/}   $D(+\infty)=\CC$, $\overline D(+\infty)$ --- расширенная комплексная плоскость, 
$\partial \overline D(+\infty)$ --- <<бесконечно удалённая точка>>.

Для $R\in \overline \RR^+$ и функции $v\colon D(R)\to \overline \RR$ 
  \begin{equation}\label{{MC}M}
{\sf M}_v(r):=\sup\bigl\{v(re^{i\theta})\bigm| \theta \in [0,2\pi)\bigr\} , \quad 0\leq r<R,
\end{equation}
--- {\it максимальная характеристика роста функции $v$  на окружностях  $\partial \overline D(r)$,} для которой  
очевидны свойства 
${\sf M}_{v^+}=({\sf M}_{v})^+=:{\sf M}_{v}^+$, ${\sf M}_{-v}^+={\sf M}_{v^-}$, ${\sf M}_{|v|}= \max\{{\sf M}_v^+,{\sf M}_{-v}^+\}$.
Будут встречаться случаи, когда для некоторых $r\in [0,R)$ значение  $v(re^{i\theta})$ определено не для всех $\theta \in [0,2\pi)$. Для таких $r$ будем считать, что  и  соответствующее значение ${\sf M}_v(r)$ не определено. 
Через  
\begin{equation}\label{{MC}C}
{\sf C}_v(r):=\frac{1}{2\pi}\int_0^{2\pi} v(re^{i\varphi})\dd \varphi 
\end{equation}
обозначаем  {\it среднее по окружности\/ $\partial \overline{D}(r)$ функции\/} $v$ при условии, что интеграл Лебега\,--\,Стилтьеса  в правой части  корректно определён значениями в $\overline \RR$.

Для $ 0\leq r\leq R\in \RR^+ $ и   меры Бореля  $\mu$ на $\overline D(R)$ возрастающая функция 
\begin{subequations}\label{murad}
\begin{align}
\mu^{\rad}(r)&:=\mu^{\rad}\bigl(\overline D(r)\bigr)
\text{ при $r\in [0, R]$}
\tag{\ref{murad}$\mu$}\label{{murad}m}
 \\
\intertext{--- {\it  радиальная считающая функция меры $\mu$}, а } 
{\sf N}_{\mu}(r,R)&\overset{\eqref{{murad}m}}{:=}\int_{r}^{R}\frac{\mu^{\rad}(t)}{t}\dd t\in \overline \RR^+
\tag{\ref{murad}N}\label{{murad}N}
\end{align}
\end{subequations} 
--- {\it  разностная  усреднённая,\/} или  {\it проинтегрированная,  радиальная считающая функция меры $\mu$}
от двух переменных $r<R$. 

Пусть  $U=u-v$ --- разность пары  субгармонических функций 
 $u\not\equiv -\infty$ и $v\not\equiv -\infty$ в окрестности круга $\overline D(R)$ с {\it мерами Рисса\/} 
соответственно 
\begin{equation}\label{df:cm}
\varDelta_u:= \frac{1}{2\pi} {\bigtriangleup}  u\geq 0, 
\quad \text{где \it $\bigtriangleup$ --- оператор Лапласа,} 
\end{equation}
и $\varDelta_v\geq 0$, являющимися {\it мерами Радона\/} в окрестности круга $\overline D(R)$.   Другими словами,   $U\not\equiv\pm\infty$ ---  {\it $\delta$-суб\-га\-р\-м\-о\-н\-и\-ч\-е\-с\-к\-ая 
 функция\/} 
с  {\it зарядом  Рисса\/} 
\begin{equation}\label{df:cmU} 
\varDelta_U:= \frac{1}{2\pi} {\bigtriangleup}  U=\varDelta_u-\varDelta_v.
\end{equation}
Различные  эквивалентные формы определения таких функций и их основные свойства приводятся и исследуются в \cite{Arsove53}, \cite{Arsove53p}, \cite[2.8.2]{Azarin}, \cite{Gr}, \cite[3.1]{KhaRoz18}. {\it Разностная   характеристика Неванлинны\/} такой   функции $U$ использовалась в нашей статье \cite{Kha20} и может быть определена 
как  функция двух переменных 
\begin{equation}\label{T}
{\sf T}_U(r,R)={\sf C}_{U^+}(R)-{\sf C}_{U^+}(r)+
{\sf N}_{\varDelta_U^-}(r,R),
\quad 0<r< R\in \RR^+, 
\end{equation}
где положительная мера $\varDelta_U^-:=\sup\{\varDelta_v, \varDelta_u\}-\varDelta_u\geq 0$ --- это {\it нижняя вариация\/} заряда Рисса  $\varDelta_{U}=\varDelta_u-\varDelta_v$ функции $U$.

Для мероморфной функции $f\neq 0,\infty$ на $\CC$
её логарифм модуля $\ln |f|\not\equiv \pm \infty$ --- нетривиальная  $\delta$-субгармоническая функция на  ${\mathbb{C}}$ и все результаты данной статьи о $\delta$-субгармонических функциях новые и для мероморфных функций с переходом, если  необходимо, к  традиционным  обозначениям 
\begin{subequations}\label{cs}
\begin{align}
\ln M(r, f)&\overset{\eqref{{TN}M},\eqref{{MC}M}}{=}{\sf M}_{\ln|f|}(r), \quad r\in {\mathbb{R}}^+,
\tag{\ref{cs}M}\label{{cs}M}\\  
m(r, f)&\overset{\eqref{{TN}m},\eqref{{MC}C}}{=}{\sf C}_{\ln^+|f|}(r),\quad r\in {\mathbb{R}}^+,
\tag{\ref{cs}m}\label{{cs}m}\\
N(R, f)-N(r, f)&\overset{\eqref{{TN}N},\eqref{{murad}N}}{=}
{\sf N}_{\varDelta_{\ln|f|}^-}(r.R)
, \quad 0<r< R\in \RR^+, 
\tag{\ref{cs}N}\label{{cs}N}\\
T(R, f)-T(r, f)&\overset{\eqref{{TN}T},\eqref{T}}{=}{\sf T}_{\ln|f|}(r,R).
\quad 0<r<R\in {\mathbb{R}}^+.
\tag{\ref{cs}T}\label{{cs}T}
\end{align}
\end{subequations}

В дальнейшем будет удобнее  использовать иную форму {\it разностной характеристики Невалинны,\/} которую  можно определить через  ${\sf T}_{U}$ из \eqref{T} в виде 
\begin{subequations}\label{rT}
\begin{align}
{\boldsymbol T}_U(r,R)&:={\sf T}_{U}(r,R)+{\sf C}_{U^+}(r)
\tag{\ref{rT}T}\label{{rT}T}
\\
&\overset{\eqref{T}}{=}{\sf C}_{U^+}(R)+{\sf N}_{\varDelta_U^-}(r,R),
\quad 0<r<R\in \RR^+,
\tag{\ref{rT}N}\label{{rT}N}
\\
\intertext{где правая часть позволяет определить и}
{\boldsymbol  T}_U(R)&:={\boldsymbol  T}_U(0,R)\overset{\eqref{{murad}N}}{:=}{\sf C}_{U^+}(R)+{\sf N}_{\varDelta_U^-}(0,R)\in \overline \RR^+.
\tag{\ref{rT}$_0$}\label{{rT}o}
\end{align}
\end{subequations}
В этом случае \eqref{{cs}T} придётся заменить  на 
\begin{equation}\label{csTT}
T(R, f)-T(r, f)+m(r,f)\overset{\eqref{{rT}T},\eqref{{cs}m}}{=}{\boldsymbol  T}_{\ln|f|}(r,R).
\end{equation}

\section{Формулировка основного результата}\label{S2}

\begin{definition}\label{DefhR}
{\it Возрастающей функции\/}   $m\colon [0,r]\to  \RR$ {\it полной вариации\/}
\begin{equation}\label{{hR}wm}
 {\tt M} :=m(r)-m(0) \in \RR^+
\end{equation}
с  {\it модулем непрерывности\/} $\omega_m\colon \RR^+\to \RR^+$, который  можно задать  как 
\begin{equation}\label{{hR}h}
\omega_m(t)\underset{t\in \RR^+}{=}\sup\bigl\{ m(x)-m(x')\bigm|x-x'\leq t, \, 0\leq x'\leq x\leq r  \bigr\}\overset{\eqref{{hR}wm}}{\subset} [0,{\tt M}],
\end{equation}
будем сопоставлять {\it диаметр стабилизации\/}
\begin{equation}\label{{hR}R}
{\sf d}_m:=\inf\bigl\{t\in \RR^+\bigm| {\omega}_{m}(t)= {\tt M}\bigr\}=\inf {\omega}_{m}^{-1}({\tt M})\leq r,
\end{equation}
где последние включение $\subset [0,{\tt M}]$ в \eqref{{hR}h} и неравенство $\leq r$ в \eqref{{hR}R} очевидны. 
Саму функция $m$ на $[0,r]$ с тем же обозначением $m$ часто будем рассматривать как {\it  возрастающую продолженную на $\RR$} постоянными  значениями $m(r)$ на луче $(r,+\infty)$ и постоянными значениями $m(0)$ на отрицательном луче $-\RR^+$, очевидно, {\it без увеличения полной вариации\/} ${\tt M}$, и тогда 
\begin{equation}\label{{hRR}h}
\omega_m(t)\overset{\eqref{{hR}h}}{\underset{t\in \RR^+}{=}}\sup_{x\in \RR}\bigl( m(x+t/2)-m(x-t/2)\bigr)\overset{\eqref{{hR}wm}}{\subset} [0,{\tt M}].
\end{equation}
{\it Носителем  непостоянства\/} $\supp m'\subset [0,r]$  продолженной таким образом функции $m$ будем называть множество, состоящее из всех точек на $\RR$, в каждой окрестности которых продолженная функция  $m$ принимает хотя бы два различных значения. Так, если продолженная  функция $m$ дифференцируема, то, очевидно, её носитель непостоянства --- это носитель ей производной. 
\end{definition}

Все {\it интегралы\/} (Римана\,-- или  Лебега\,--\,){\it Стилтьеса\/}
с нижним пределом интегрирования    $a\in \RR$ и и верхним пределом интегрирования $b\geq a$  понимаем как интеграл по отрезку $[a,b]$, если не оговорено иное.

\begin{theorem}\label{th1} Пусть\/  $0< r<R\in \RR^+$, 
$m\colon [0,r]\to\RR$ --- возрастающая функция.   Если  для модуля непрерывности $\omega_m$ выполнено условие Дини
\begin{equation}\label{{hR}i}
\int_0^{4R}\frac{{\omega}_{m}(t)}{t}\dd t<+\infty,
\end{equation}
то  для любой $\delta$-субгармонической функции $U\not\equiv\pm\infty$ в окрестности замкнутого круга 
$\overline D(R)$ существует интеграл Лебега\,--\,Стилтьеса   с верхней оценкой 
\begin{equation}\label{U}
\int_0^r {\sf M}_{U}^+(t)\dd m(t)\leq \frac{6R}{R-r}
{\boldsymbol   T}_U(r,R) \max\biggl\{{\tt M}, \int_0^{{\sf d}_m}\ln \frac{4R}{t}\dd {\omega}_{m}(t)\biggr\}, 
\end{equation}
 где первый аргумент $r$ в ${\boldsymbol   T}_U(r,R)$ 
можно заменить на любое  $r_0\in [0,r]$, а   последний  интеграл Римана\,--\,Стилтьеса в \eqref{U} под операцией $\max$ --- на сумму
\begin{equation}\label{kint}
\int_0^{{\sf d}_m} \frac{{\omega}_{m}(t)}{t}\dd t+
{\tt M}\ln \frac{4R}{{\sf d}_m}\geq \int_0^{{\sf d}_m} \ln\frac{4R}{t}\dd {\omega}_{m}(t).
\end{equation} 
\end{theorem}
\begin{remark}\label{rem1}
Конечно, в классическом  условии Дини в верхнем  пределе интеграла, в отличие от \eqref{{hR}i},  ставится сколь угодно малое  строго положительное число, но для ссылок будет удобна именно такая эквивалентная форма. Фактически единственное в основной теореме условие типа Дини обязано быть выполненным, как показывает следующий контрпример.

\begin{example} Пусть $r=2$ и возрастающая функция $m\colon [0,2]\to \RR$ такова, что 
\begin{equation}\label{m1}
\int_0^1\frac{m(1+t)-m(1-t)}{t}\dd t=+\infty,
\end{equation}
что означает нарушение условия Дини для функции $m$ в точке $1$  в самой  общей форме без каких-либо особых дополнительных требований. Рассмотрим мероморфную функция $z\longmapsto 5/(z-1)$ на $\CC$. Тогда для  {\it супергармонической,\/} а значит, и $\delta$-субгармонической  функции  
\begin{equation}\label{U5}
U(z)=\ln \Bigl|\frac{5}{z-1}\Bigr|=\ln \frac{5}{|z-1|}, \quad z\in \CC,
\end{equation}
положительной на $\overline D(4)$, имеем  
\begin{multline*}
\int_0^2{\sf M}_{U}^+(x)\dd m(x)=\int_0^2\ln \frac{5}{|x-1|}\dd m(x)
=\int_0^1\ln \frac{5}{t}\dd \,\bigl(m(1+t)-m(1-t)\bigr)\\
=\bigl(m(2)-m(0)\bigr)\ln 5+\int_0^1\frac{m(1+t)-m(1-t)}{t}\dd t\overset{\eqref{m1}}{=}+\infty,
\end{multline*}
вследствие  чего конечные оценки сверху для  интеграла слева   невозможны.
\end{example} 

Этот же контрпример указывает, что  интегрирование 
максимальной характеристики роста произвольной $\delta$-субгармонической функции имеет смысл рассматривать только по возрастающей  функции $m$, в классическом разложении которой не должно быть {\it функции скачков,\/} а могут быть только {\it абсолютно непрерывная\/} и {\it сингулярная слагаемые\/} относительно меры Лебега $\mes$ на $[0,r]$. 
\end{remark}
\begin{remark}\label{rem2}
 Оценки интегралов от {\it максимальной радиальной характеристики роста\/} для общих $\delta$-субгармонических функций от более чем двух вещественных переменных или для общих мероморфных  функций или  разностей плюрисубгармонических более чем одной комплексной переменной невозможны. Приведём  примеры таких невырожденных функций с тождественно равной $+\infty$  максимальной характеристикой роста.   
\begin{example}
Пусть $n\geq 3$, $x=(x_1,\dots, x_n)\in \RR^n$. 

Для  супергармонической на $\RR^3$  функции  
$$
u(x)=u(x_1,x_2, x_3):=\begin{cases}
-\ln\sqrt{x_1^2+x_2^2} &\text{при $x_1^2+x_2^2\neq 0$}\\
+\infty  &\text{при $x_1^2+x_2^2= 0$}
\end{cases}
$$
и для супергармонической на $\RR^n$ при $n\geq 4$ функции  
$$
u(x)=u(x_1, \dots, x_n):=\begin{cases}
1\biggm/\Biggl(\sqrt{\sum\limits_{j=1}^{n-1}x_j^2}\Biggr)^{n-3}&\text{при $\sum\limits_{j=1}^{n-1}x_j^2\neq 0$}\\
+\infty  &\text{при $\sum\limits_{j=1}^{n-1}x_j^2= 0$}
\end{cases}
$$
максимум её на сфере с центром в нуле  любого радиуса $r\in \RR^+$ равен $+\infty$.  
\end{example}
\begin{example}
Пусть $n\geq 2$, $z=(z_1,\dots, z_n)\in \CC^n$. Для мероморфной на $\CC^n$ функции $f(z)=f(z_1, \dots, z_n)={1}/{z_1}$ логарифм максимума её модуля  на каждой сфере радиуса $r$ с центром в нуле равен $+\infty$.  
\end{example}

\end{remark}

\section{Разностные характеристики Неванлинны}
Две $\delta$-субгармонические функции $U=u-v\not\equiv \pm \infty$ и $U_1=u_1-v_1\not\equiv \pm \infty$, представленные разностями пар субгармонических функций $u,v,u_1,v_1\not\equiv -\infty$ на открытом множестве, {\it равны\/} на этом множестве, если $u+v_1=u_1+v$ на нём \cite{Arsove53}, \cite{Arsove53p}, \cite[2.8.2]{Azarin}, \cite{Gr}, \cite[3.1]{KhaRoz18}. Исходное определение разностной характеристики Неванлинны ${\sf T}_U$ в  \cite{Kha20} вводится иначе.  Для $\delta$-суб\-г\-а\-р\-м\-о\-н\-и\-ч\-е\-с\-к\-ой функции $U\not\equiv \pm\infty$ на круге $\overline D(R)$, т.е. в некоторой окрестности круга  $\overline D(R)$, с зарядом Рисса $\varDelta_U$ существуют {\it канонические   представления\/} $U=u_U-v_U$, где $u_U\not\equiv -\infty$  и $v_U\not\equiv -\infty$  --- субгармонические функции на $\overline D(R)$ с мерами Рисса соответственно 
$\varDelta_{u_U}=\varDelta_U^+$ и $\varDelta_{v_U}=\varDelta_U^-$. Канонические представления  определены с точностью до общего гармонического слагаемого. Очевидны равенства 
\begin{equation}\label{uvU}
\sup\{u_U,v_U\}=\sup\{u_U-v_U, 0\}+v_U=U^++v_U,
\end{equation}
где слева --- субгармоническая функция. Если  $v\not\equiv -\infty$ --- субгармоническая функция на $\overline D(R)$, то по  классической формуле Пуассона\,--\,Йенсена\,--\,При\-в\-а\-л\-о\-ва для круга \cite[4.5]{Rans}, \cite[4]{Pri35}, \cite[гл. II, \S~2]{Privalov}
\begin{equation}\label{CN}
{\sf C}_v(R)-{\sf C}_v(r)
={\sf N}_{\varDelta_v}(r,R)
\quad\text{для всех  $0<r<R<+\infty$}.
\end{equation}

 В  \cite[(15)]{Kha20} определена разностная характеристика  Неванлинны 
\begin{subequations}\label{TO}
\begin{align}
{\sf T}_U(r,R)&:={\sf C}_{\sup\{u_U,v_U\}}(R)-{\sf C}_{\sup\{u_U,v_U\}}(r)
\tag{\ref{TO}T}\label{{TO}T}\\
&\overset{\eqref{uvU}}{=}{\sf C}_{U^+}(R)-{\sf C}_{U^+}(r)+{\sf C}_{v_U}(R)-{\sf C}_{v_U}(r)
\tag{\ref{TO}C}\label{{TO}C}
\\
&\overset{\eqref{CN}}{=}{\sf C}_{U^+}(R)-{\sf C}_{U^+}(r)+{\sf N}_{\varDelta_U^-}(r,R), \quad 0<r\leq R\in \RR^+. 
\tag{\ref{TO}N}\label{{TO}N}
\end{align}
\end{subequations}
В заключительной части \cite[формула (15)]{Kha20}, где  это определение  преобразуется к виду \eqref{T}, к сожалению,    
была  допущена описка: вместо необходимого, как в правой части  \eqref{T}, слагаемого ${\sf N}_{\varDelta_U^-}(r,R)$ с нижней вариацией $\varDelta_U^-$  ошибочно написано  ${\sf N}_{\varDelta_U^+}(r,R)$ с   верхней вариацией $\varDelta_U^+$. Впрочем, на правильности выводов  результатов из \cite{Kha20}  это не сказалось, поскольку всюду в \cite{Kha20} в доказательствах применяется  корректно записанное   определение.

Cвойства разностной характеристики Неванлинны ${\sf T}_U$, которые сразу следуют из свойств среднего по окружности \eqref{{MC}C} \cite[теорема 2.6.8]{Rans},   объединим в 
\begin{propos}\label{pr1} Пусть $0<R\in \RR^+$, а   $U\not\equiv\pm\infty$ ---  $\delta$-суб\-га\-р\-м\-о\-н\-и\-ч\-е\-с\-к\-ая   функция  на $\overline D(R)$. Разностная характеристика Неванлинны ${\sf T}_U$ из \eqref{T}
\begin{enumerate}[{\rm (i)}]
\item\label{T1}   не зависит от её представления в виде разности субгармонических,  
\item\label{T0} положительна и непрерывна по каждой переменной $0<r<R<+\infty$,
\item\label{T2} возрастает и выпукла относительно логарифма $\ln$ по второму аргументу,  а также 
 убывает и вогнута относительно логарифма $\ln$ по  первому; 
\item\label{T3} существует конечный или бесконечный предел
\begin{equation}\label{sfT0}
\lim_{0<r\to 0} {\sf T}_U(r,R)=:{\sf T}_U(0,R)=:{\sf T}_U(R)\in \overline \RR^+,
\end{equation}
\item\label{T4} ${\sf T}_{U}={\sf T}_{-U}$.
\item\label{T5} если функция $u$ субгармоническая на $\overline D(R)$ с мерой Рисса $\varDelta_u$, то 
\begin{equation*}
{\sf T}_{u}(r,R)={\sf C}_{u^+}(R)-{\sf C}_{u^+}(r)={\sf C}_{u^-}(R)-{\sf C}_{u^-}(r)+
{\sf N}_{\varDelta_u}(r,R)={\sf T}_{-u}(r,R).
\end{equation*} 
\item\label{T6} если $U=u-v$, а для  субгармонических функций $u$ и $v$ имеем  $u(0)\neq -\infty$ или $v(0)\neq -\infty$, то предел в \eqref{sfT0} конечный и  ${\sf T}_U(R)\in \RR^+$.

\end{enumerate}
\end{propos}
\begin{proof} Свойство \eqref{T1} следует из \eqref{{TO}N}; \eqref{T0} и \eqref{T2}  --- из 
\eqref{{TO}T} в силу возрастания, непрерывности  и выпуклости относительно $\ln$ среднего по окружности 
${\sf C}_{\sup\{u_U,v_U\}}$ субгармонической функции $\sup\{u_U,v_U\}$; свойство \eqref{T3}
сразу следует из свойства \eqref{T2};  \eqref{T4} --- из возможности переставить 
$u_U$ и $v_U$ в \eqref{{TO}T}; \eqref{T5} --- из \eqref{T4} и определений \eqref{TO} или \eqref{CN} для промежуточного равенства. В условиях п.~\eqref{T6} для канонического представления  $U=u_U-v_U$ тем более $u_U(0)\neq -\infty$ или 
$v_U\neq -\infty$. Отсюда $\sup\{u_U,v_U\}(0)\neq -\infty$ и для субгармонической функции $\sup\{u_U,v_U\}$ существует 
предел \cite[теорема 2.6.8]{Rans}
$$
\lim_{0<r\to 0}{\sf C}_{\sup\{u_U,v_U\}}(r)=\sup\{u_U,v_U\}(0)\in \RR,
$$
откуда по определению \eqref{{TO}T}  предел в \eqref{sfT0} конечный. 
\end{proof}
Свойства характеристики ${\boldsymbol   T}_U$ из \eqref{rT} объединим в 
\begin{propos}\label{pr2} 
В  условиях предложения\/ {\rm \ref{pr1}} разностная характеристика Неванлинны вида ${\boldsymbol   T}_U$ из \eqref{rT}
\begin{enumerate}[{\rm (i)}]
\item\label{bT1}   не зависит от её представления в виде разности субгармонических,  
\item\label{bT0} положительна и непрерывна по каждой переменной	$0<r<R<+\infty$,
\item\label{bT2} возрастает и выпукла относительно  $\ln$ по второму аргументу,  а также 
 убывает и вогнута относительно $\ln$ по  первому, в частности, 
\begin{equation}\label{{TT}r} 
{\boldsymbol  T}_U(r,R)\leq {\boldsymbol  T}_U(r_0,R)
\text{ при всех  $0\leq r_0\leq r<R\in \RR^+$};
\end{equation}

\item\label{bT4} ${\boldsymbol  T}_{U}(r,R)={\boldsymbol T}_{-U}(r,R)+{\sf C}_U(r)$ при $0<r<R$;

\item\label{bT5} если функция $u$ субгармоническая на $\overline D(R)$ с мерой Рисса $\varDelta_u$, то 
\begin{multline*}
{\boldsymbol T}_{u}(r,R)\overset{\eqref{{rT}N}}{=}{\sf C}_{u^+}(R)={\sf C}_{u}(R)+{\sf C}_{u^-}(R)
\\={\sf C}_{(-u)^+}(R)+{\sf N}_{\varDelta_u}(r,R)+{\sf C}_{u}(r)\overset{\eqref{{rT}N}}{=}{\boldsymbol T}_{-u}(r,R)+{\sf C}_{u}(r).
\end{multline*} 

\item\label{bT6} 
 если  в некотором  представлении $U=u-v$  разностью  субгармонических функций  $u$ и $v$ имеем $v(0)\neq -\infty$, то   существуют  конечные значения 
\begin{equation}\label{Ntl}
\lim_{0<r\to 0}{\sf N}_{\varDelta_U^-}(r,R)\in \RR^+, \quad  {\boldsymbol T}(R)\overset{\eqref{{rT}o}}{=}{\boldsymbol T}(0,R)\in \RR^+ .
\end{equation}
\end{enumerate}
\end{propos}
\begin{proof} Свойства \eqref{bT1} и  \eqref{bT0} пока без непрерывности, а также часть \eqref{bT2} по второму аргументу следуют из  определения \eqref{{rT}T} и пп. \eqref{T1}--\eqref{T2}  предложения \ref{pr1}. Часть \eqref{bT2} по первому аргументу вместе с \eqref{{TT}r} следует из представления \eqref{{rT}N} и соответствующих свойств 
разностной   усреднённой проинтегрированной  радиальной считающей функции  меры. 
Непрерывность в \eqref{bT0} --- следствие выпуклости относительно $\ln$ по $R$ и вогнутости относительно $\ln$ по $r$.
Используя  п.~\eqref{T4}
предложения \ref{pr1} и  определение  \eqref{{rT}T}, из равенств 
\begin{multline*}
{\boldsymbol T}_U(r,R)\overset{ \eqref{{rT}T}}{=}{\sf T}_{U}(r,R)+{\sf C}_{U^+}(r)={\sf T}_{-U}(r,R)+{\sf C}_{U^+}(r)
\\
\overset{ \eqref{{rT}T}}{=}{\boldsymbol T}_{-U}(r,R)-{\sf C}_{(-U)^+}(r)+{\sf C}_{U^+}(r) 
={\boldsymbol T}_{-U}(r,R)+{\sf C}_{U^+}(r) -{\sf C}_{U^-}(r),
\end{multline*}
получаем в точности п.~\eqref{bT4}  предложения \ref{pr2}.
Из п.~\eqref{bT4} и \eqref{{rT}N}  следует п.~\eqref{bT5}.  

Если $v(0)\in \RR$, то существует предел
$\lim\limits_{0<r\to 0}{\sf C}_{v}(r)=v(0)\in \RR$  \cite[теорема 2.6.8]{Rans}, 
откуда согласно \eqref{CN} существует конечный предел 
$\lim\limits_{0<r\to 0}{\sf N}_{v}(r,R)\in \RR^+$.
Но $\varDelta_v\geq \varDelta_U^-$, следовательно, конечен  и предел в \eqref{Ntl}, откуда
${\boldsymbol T}(0, R)\overset{\eqref{{rT}o}}\in \RR^+$. 
\end{proof}

\section{Элементарная оценка $\delta$-субгармонической функции в круге} 
\begin{propos}\label{PJDR} Пусть $0<r<R\in \RR^+$,  $U\not\equiv\pm\infty$ --- 
 $\delta$-суб\-г\-а\-р\-м\-о\-н\-и\-ч\-е\-с\-к\-ая  функция на  круге $\overline D(R)\subset \CC$ радиуса $R>0$ с зарядом Рисса $\varDelta_U$. Тогда 
\begin{equation}\label{estU}
U(w)\leq \frac{R+r}{R-r}{\sf C}_{U^+}(R) 
+\int_{D(R)}\ln \frac{2R}{|w-z|}\dd \varDelta_{U}^-(z) 
\quad\text{при всех $w\in \overline D(r)$.}
\end{equation}
\end{propos}
\begin{proof}
Дважды применяя формулу Пуассона\,--\,Йенсена \cite[4.5]{Rans} для круга $D(R)$ к паре 
субгармонических функций $u\not\equiv -\infty$ и $v\not\equiv -\infty$ с мерами Рисса соответственно $\varDelta_u$ и $\varDelta_v$, представляющих $U=u-v$  на круге $\overline D(R)$,    после вычитания имеем равенство    
\begin{equation*}
U(w)=\frac{1}{2\pi}\int_0^{2\pi}U(Re^{i\varphi})\Re \frac{Re^{i\varphi}+w}{Re^{i\varphi}-w}\dd \varphi -\int_{D(R)}\ln \Bigl|\frac{R^2-z\bar w}{R(w-z)}\Bigr|\dd\, (\varDelta_u-\varDelta_v)(z) 
\end{equation*}
при всех $w\in \overline D(r)$ с положительным ядром Пуассона
\begin{equation*}
0\leq \Re \frac{Re^{i\varphi}+w}{Re^{i\varphi}-w}\leq 
\Bigl|\frac{Re^{i\varphi}+w}{Re^{i\varphi}-w}\Bigr|\leq
\frac{R+|w|}{R-|w|}\leq \frac{R+r}{R-r}\text{ при $|w|\leq r$}
\end{equation*}
и положительной функцией Грина
\begin{equation*}
0\leq \ln \Bigl|\frac{R^2-z\bar w}{R(w-z)}\Bigr|\leq
\ln \frac{R^2+|z\bar w|}{R|w-z|}\leq
\ln \frac{R^2+R^2}{R|w-z|}= \ln \frac{2R}{|w-z|}
\text{ при $z\in  D(R)$}.
\end{equation*}
Отсюда для заряда Рисса $\varDelta_{U}=\varDelta_u-\varDelta_v$ функции $U$ с положительной вариацией 
$\varDelta_{U}^+=(\varDelta_u-\varDelta_v)^+$ и отрицательной вариацией $\varDelta_{U}^-=(\varDelta_u-\varDelta_v)^-$
получаем 
\begin{multline*}
U(w)\leq \frac{R+r}{R-r}{\sf C}_{U^+}(R) 
+\int_{D(R)}\ln \Bigl|\frac{R^2-z\bar w}{R(w-z)}\Bigr|\dd (\varDelta_{U}^--\varDelta_{U}^+)(z) 
\\
\leq  \frac{R+r}{R-r}{\sf C}_{U^+}(R)  +\int_{D(R)}\ln \frac{2R}{|w-z|}\dd \varDelta_{U}^-(z) 
\quad\text{при всех $w\in \overline D(r)$}.
\end{multline*}
\end{proof}

\section{Доказательство  теоремы \ref{th1}}
Конкретный выбор промежуточного $R_*$ между $r$ и $R$ укажем позже, а пока 
\begin{equation}\label{rR}
0<r<R_*<R, \quad \overline D(r)\subset D(R_*)\subset D(R)\subset \overline D(R). 
\end{equation} 
Рассмотрим сначала случай, когда функция $U$ представима в виде разности $U:=u-v$ пары бесконечно дифференцируемых субгармонических функций  $u$ и $v$ в окрестности круга $\overline D(R)$ 
с мерами Рисса $\varDelta_u$ и $\varDelta_v$ с гладкими плотностями относительно плоской меры Лебега. Это позволит пока не утруждаться обоснованием выкладок. По предложению \ref{PJDR} из неравенства \eqref{estU} для  $R_*$ вместо  $R$ при всех $w:=te^{i\theta}\in \in \overline D(r)$ получаем 
\begin{multline*}
U(te^{i\theta})\leq \frac{R_*+r}{R_*-r}{\sf C}_{U^+}(R_*) 
+\int_{D(R_*)}\ln \frac{2R_*}{|te^{i\theta}-z|}\dd \varDelta_{U}^-(z) \\
\leq  \frac{R_*+r}{R_*-r}{\sf C}_{U^+}(R_*)   +\int_{D(R_*)}\ln \frac{2R_*}{\bigl|t-|z|\bigr|}\dd \varDelta_{U}^-(z) 
\\=\frac{R_*+r}{R_*-r}{\sf C}_{U^+}(R_*) +\int_0^{R_*}\ln \frac{2R_*}{|t-x|}\dd \,(\varDelta_U^-)^{\rad}(x),
\end{multline*}
где использовано обозначение  \eqref{{murad}m} для радиальной считающей функции $(\varDelta_U^-)^{\rad}$ меры $\varDelta_U^-$, а  слагаемые в правой части неравенства положительны и не зависят от $\theta\in[0,2\pi)$. Поэтому переход  к точной верхней грани по $\theta\in[0,2\pi)$  в левой части полученного неравенства, а затем  к положительной части  влечёт за собой  неравенство
\begin{equation*}
{\sf M}_{U}^+(t)\leq \frac{R_*+r}{R_*-r}{\sf C}_{U^+}(R_*)+
\int_0^{R_*}\ln \frac{2R_*}{\bigl|t-x\bigr|}\dd \,(\varDelta_U^-)^{\rad}(x)\quad\text{при всех $t\in[0, r]$.}
\end{equation*}
Интегрируя по возрастающей функции $m$ и используя теорему Фубини о повторных интегралах, приходим к соотношениям  
\begin{multline}\label{Iln}
\int_0^r{\sf M}_{U}^+(t)\dd m(t)\\
\leq \int_0^r\frac{R_*+r}{R_*-r} 
{\sf C}_{U^+}(R_*)\dd m(t)
+\int_0^r\int_0^{R_*}\ln \frac{2R_*}{|t-x|}\dd \,(\varDelta_U^-)^{\rad}(x)\dd m(t)\\
=\frac{R_*+r}{R_*-r} {\sf C}_{U^+}(R_*)\int_0^r\dd m(t)
+ \int_0^{R_*}\int_0^r\ln \frac{2R_*}{|t-x|}\dd m(t)\dd \,(\varDelta_U^-)^{\rad}(x)
\\
\overset{\eqref{{hR}wm}}{\leq} \frac{R_*+r}{R_*-r} {\sf C}_{U^+}(R_*){\tt M}+(\varDelta_U^-)^{\rad}(R_*)
\sup_{x\in [0,R_*]}\int_0^r\ln \frac{2R_*}{|t-x|}\dd m(t)\\
\overset{\eqref{rR}}{\leq} \frac{R_*+r}{R_*-r} {\sf C}_{U^+}(R_*){\tt M}+(\varDelta_U^-)^{\rad}(R_*)
\sup_{x\in [0,R]}\int_0^r\ln \frac{2R}{|t-x|}\dd m(t),
\end{multline}
где при последнем переходе использован выбор $r<R_*\overset{\eqref{rR}}{<}R$.
При замене $s:=|t-x|$ в последнем   интеграле в \eqref{Iln}, используя 
продолженную на $\RR$ без изменения полной вариации  возрастающую функции $m$, как в определении \ref{DefhR},  
получаем представление
\begin{multline*}
\int_0^r\ln \frac{2R}{|t-x|}\dd m(t)=\int_0^{R}\ln \frac{2R}{s}\dd \,\bigl(m(x+s)
-m(x-s)\bigr)\\
=\int_0^{2R}\ln \frac{4R}{t}\dd \,\bigl(m(x+t/2)
-m(x-t/2)\bigr)
=\int_0^{4R}\ln \frac{4R}{t}\dd \,\bigl(m(x+t/2)
-m(x-t/2)\bigr),
\end{multline*}
где последнее равенство следует из постоянства при указанном продолжении функции 
$m(x+t/2)-m(x-t/2)$ для  каждого фиксированного  $x\in [0,R]$ при всех $t\geq 2R$.  Интегрирование  по частям последнего  интеграла 
даёт  равенство 
\begin{multline*}
\int_0^r\ln \frac{2R}{|t-x|}\dd m(t)=
-\lim_{0<t\to 0}\bigl(m(x+t/2)
-m(x-t/2)\bigr)\ln \frac{4R}{t}\\
+ \int_0^{4R}\frac{m(x+t/2)
-m(x-t/2)}{t}\dd t,
\end{multline*}
Здесь,   по определению модуля непрерывности  \eqref{{hR}h} и условию Дини  \eqref{{hR}i} имеем 
\begin{equation*}
\int_0^{4R}\frac{m(x+t/2)
-m(x-t/2)}{t}\dd t\overset{\eqref{{hRR}h}}{\leq} \int_0^{4R}\frac{{\omega}_{m}(t)}{t}\dd t\overset{\eqref{{hR}i}}{<}+\infty,
\end{equation*}
 откуда  последний предел   $-\lim\limits_{0<t\to 0}\dots$ равен нулю \cite[предложение 2.2]{KhI}, и  
\begin{multline*}
\sup_{x\in [0,R]}\int_0^r\ln \frac{2R}{|t-x|}\dd m(t)=
\sup_{x\in [0,R]} \int_0^{4R}\frac{m(x+t/2)
-m(x-t/2)}{t}\dd t\\
\leq \int_0^{4R}\sup_{x\in \RR} \bigl(m(x+t/2)
-m(x-t/2)\bigr)\frac{\dd t}{t}
\overset{\eqref{{hRR}h}}{=}\int_0^{4R}\frac{{\omega}_{m}(t)}{t}\dd t.
\end{multline*}
Таким образом,  неравенство \eqref{Iln} переходит в неравенство
\begin{equation}\label{iMd-}
\int_0^r{\sf M}_{U}^+(t)\dd m(t)
\leq \frac{R_*+r}{R_*-r} {\sf C}_{U^+}(R_*){\tt M}+(\varDelta_U^-)^{\rad}(R_*)
\int_0^{4R}\frac{{\omega}_{m}(t)}{t}\dd t. 
\end{equation}
\begin{lemma}\label{lemin} Для любой $\delta$-субгармонической функции $U\not\equiv \pm\infty$ на $\overline D(R)$ 
выполнено неравенство \eqref{iMd-}  с конечной правой частью. 
\end{lemma}
\begin{proof}[леммы \ref{lemin}]
Перенесём теперь последнее неравенство на случай {\it произвольной  $\delta$-субгармонической функции $U=u-v$,}
представленной в виде разности субгармонических функций $u\not\equiv-\infty$ и $v\not\equiv -\infty$ на $\overline D(R)$. 
По \cite[теорема 2.7.2]{Rans} существуют убывающие последовательности бесконечно дифференцируемых субгармонических функций $(u_n)_{n\in \NN}$ и $(v_k)_{k\in \NN}$ на $\overline D(R)$, стремящиеся поточечно к $u$ и $v$
соответственно, а  их меры Рисса $\varDelta_{u_n}$ и $\varDelta_{v_k}$ при этом $*$-слабо сходятся соответственно к мерам Рисса   $\varDelta_{u}$ и $\varDelta_{v}$ на $D(R)$ \cite[A4]{Rans}.   Для каждой бесконечно дифференцируемой $\delta$-субгармонической функции $U^n_k:=u_n-v_k$ на $D(R)$   функция ${\sf M}_{U^n_k}^+$ непрерывна  на $[0,R]$.
При фиксированном $k\in \NN$ рассмотрим $\delta$-субгармоническую  функцию $U_k$, полученную как поточечный предел
 {\it убывающей последовательности\/}  $U^n_k$ при $n\to \infty$. По определению среднего по окружности 
и теореме Лебега  о мажорируемой сходимости интегралов  ${\sf C}_{{U_k}^+}(R_*)$ --- это предел убывающей последовательности средних по окружности $\partial \overline D(R_*)$ от функций $U^n_k$. Кроме того,  убывающая последовательность {\it непрерывных функций\/} ${\sf M}_{U^n_k}^+\geq {\sf M}_{U_k}^+$  при $n\to \infty$ поточечно 
стремится к {\it полунепрерывной сверху функции\/} ${\sf M}_{U_k}^+$ на $[0,r]$. Наконец, 
из $*$-слабо сходимости зарядов  Рисса $\varDelta_{u_n}-\varDelta_{v_k}$ 
к $\varDelta_{u}-\varDelta_{v_k}$  для любого числа   $R_0\in (R^*,R)$  имеем 
$(\varDelta_{U^k}^-)^{\rad}(R_0)\geq\bigl( (\varDelta_{u} -\varDelta_{v_k})^-\bigr)^{\rad}(R_*)$ при всех $k\in \NN$. Отсюда применение  \eqref{iMd-} к каждой функции $U^n_k$ с последующим предельным переходом по $n\to \infty$ даёт для каждого $k\in \NN$ неравенство 
\begin{equation}\label{iMd-R}
\int_0^r{\sf M}_{U_k}^+(t)\dd m(t)
\leq \frac{R_*+r}{R_*-r} {\sf C}_{U_k^+}(R_*){\tt M}+(\varDelta_{U_k}^-)^{\rad}(R_0)
\int_0^{4R}\frac{{\omega}_{m}(t)}{t}\dd t 
\end{equation}
при $R_*<R_0<R$. Здесь  для {\it возрастающей последовательности\/}  $\delta$-суб\-г\-а\-р\-м\-о\-н\-и\-че\-с\-к\-их функций $U_k$
имеем ${\sf C}_{U_k^+}(R_*)\leq {\sf C}_{U^+}(R_*)<+\infty$ для всех $k\in \NN$, 
а заряды Рисса $\varDelta_{U_k}$ по-прежнему   $*$-слабо сходятся к  заряду Рисса $\varDelta_{U}$, откуда 
$(\varDelta_{U}^-)^{\rad}(R')\geq (\varDelta_{U_k}^-)^{\rad}(R_0)$ при любом $ R'\in (R_0,R)$. Таким образом, из \eqref{iMd-R} 
\begin{equation}\label{iMd-RR}
\int_0^r{\sf M}_{U_k}^+(t)\dd m(t)
\leq \frac{R_*+r}{R_*-r} {\sf C}_{U^+}(R_*){\tt M}+(\varDelta_{U}^-)^{\rad}(R')
\int_0^{4R}\frac{{\omega}_{m}(t)}{t}\dd t 
\end{equation}
при любом $k\in \NN$ и любом выборе числа $R'\in (R_*,R)$ в силу произвола в выборе числа 
$R_0\in (R_*,R)$ и  независимости левой части от $R'$. Но радиальная считающая функция меры Радона непрерывна справа, поэтому $R'$ в правой части \eqref{iMd-RR} можно заменить на $R_*$. 
Правая часть   в \eqref{iMd-RR} по условию \eqref{{hR}i}
ограничена сверху числом, не зависящим от $k\in \NN$, и, применяя теорему Леви о монотонной сходимости, получаем 
\begin{equation}\label{iMd-RR*}
\int_0^r\lim_{k\to \infty} {\sf M}_{U_k}^+(t)\dd m(t)
\leq \frac{R_*+r}{R_*-r} {\sf C}_{U^+}(R_*){\tt M}+(\varDelta_{U}^-)^{\rad}(R_*)
\int_0^{2R}\frac{{\omega}_{m}(t)}{t}\dd t ,
\end{equation}
где подынтегральный предел --- это  функция ${\sf M}_{U}^+$, поскольку в силу возрастания последовательности функций ${\sf M}_{U_k}^+$ этот предел можно рассматривать как   повторные точные верхние грани, которые можно переставлять.  

Лемма \ref{lemin} доказана. 
\end{proof}

Для последнего конечного интеграла из \eqref{iMd-}, интегрируя по частям, имеем 
\begin{multline*}
\int_0^{4R}\frac{{\omega}_{m}(t)}{t}\dd t= 
\int_0^{4R}{\omega}_{m}(t)\dd \,\Bigl(-\ln \frac{4R}{t}\Bigr)
\\= -\lim_{0<t\to 0} {\omega}_{m}(t)\ln \frac{4R}{t}
+ \int_0^{4R}\ln \frac{4R}{t}\dd {\omega}_{m}(t)
=\int_0^{4R}\ln \frac{4R}{t}\dd {\omega}_{m}(t), 
\end{multline*}
поскольку  промежуточный предел  $-\lim\limits_{0<t\to 0}\dots$ равен нулю  
ввиду конечности интеграла в левой части  \cite[предложение 2.2]{KhI}.
При этом по определению \eqref{{hR}R} диаметра стабилизации ${\sf d}_m\leq r$ модуль непрерывности  ${\omega}_{m}$ постоянен при всех $t> {\sf d}_m$. Таким образом, 
верхний предел интегрирования в последнем интеграле можно заменить на диаметр стабилизации ${\sf d}_m$
и   \eqref{iMd-} переходит в неравенство 
\begin{equation}\label{iMd}
\int_0^r{\sf M}_{U}^+(t)\dd m(t)
\leq \frac{R_*+r}{R_*-r} {\sf C}_{U^+}(R_*){\tt M}+(\varDelta_U^-)^{\rad}(R_*)
\int_0^{{\sf d}_m}\ln \frac{4R}{t}\dd {\omega}_{m}(t)<+\infty.
\end{equation}
В то же время  для меры $\varDelta_U^-$ на  замкнутом круге  $\overline D(R)$ имеем
\cite[лемма 1]{Kha20}
\begin{equation}\label{N}
(\varDelta_U^-)^{\rad}(R_*)\overset{\eqref{{murad}N}}{\leq}  \frac{R}{R-R_*}{\sf N}_{\varDelta_U^-}(R_*,R)
\overset{\eqref{{murad}N}}{\leq} \frac{R}{R-R_*}{\sf N}_{\varDelta_U^-}(r,R).
\end{equation} 
При выборе $R_*:=\frac12(R+r)$ в \eqref{rR} как среднего арифметического $r$ и $R$ из  
\begin{equation*}
\frac{R_*+r}{R_*-r}=\frac{R+3r}{R-r}\leq 4\frac{R}{R-r}, 
\quad \frac{R}{R-R_*}=\frac{2R}{R-r}
\end{equation*}
и  \eqref{iMd} в сочетании  с \eqref{N} получаем 
\begin{multline}\label{iMdR}
\int_0^r{\sf M}_{U}^+(t)\dd m(t)
\leq \frac{R}{R-r} \biggl(4{\sf C}_{U^+}(R_*){\tt M}+
2 {\sf N}_{\varDelta_U^-}(r,R)
\int_0^{{\sf d}_m}\ln \frac{4R}{t}\dd {\omega}_{m}(t)\biggr)\\
\leq 
\frac{R}{R-r} \Bigl(4{\sf C}_{U^+}(R_*)+
2 {\sf N}_{\varDelta_U^-}(r,R)\Bigr)
\max\biggl\{{\tt M}, \int_0^{{\sf d}_m}\ln \frac{4R}{t}\dd {\omega}_{m}(t)\biggr\}
<+\infty.
\end{multline} 
Здесь для круглой скобки в правой части, не зависящей от функции $m$, имеем
\begin{equation*}
4{\sf C}_{U^+}(R_*)+2 {\sf N}_{\varDelta_U^-}(r,R)
\overset{\eqref{rT}}{\leq} 4{\boldsymbol   T}_{U}(r,R_*)+2 {\boldsymbol   T}_{U}(r,R)
\leq 6{\boldsymbol   T}_{U}(r,R). 
\end{equation*}
где использовано предложение \ref{pr2}\eqref{bT2} о возрастания разностной характеристики Неванлинны \eqref{rT} по второму аргументу. Используя это неравенство для оценки правой части   \eqref{iMdR}, имеем
\begin{equation*}
\int_0^r{\sf M}_{U}^+(t)\dd m(t)\leq 
\frac{6R}{R-r} {\boldsymbol  T}_{U}(r, R)
\max\biggl\{{\tt M}, \int_0^{{\sf d}_m}\ln \frac{4R}{t}\dd {\omega}_{m}(t)\biggr\}
<+\infty,
\end{equation*} 
при всех  $r_0\in (0,r]$, что доказывает неравенство \eqref{U}. 
Для последнего интеграла интегрированием по частям   получаем 
\begin{equation*}
\int_0^{{\sf d}_m} \ln\frac{4R}{t}\dd {\omega}_{m}(t)=
\int_0^{{\sf d}_m} \frac{{\omega}_{m}(t)}{t}\dd t+
{\omega}_{m}({\sf d}_m)\ln \frac{4R}{{\sf d}_m}
\overset{\eqref{{hR}h}}{\leq}  \int_0^{{\sf d}_m} \frac{{\omega}_{m}(t)}{t}\dd t+
{ {\tt M}}\ln \frac{4R}{{\sf d}_m},
\end{equation*} 
 что даёт \eqref{kint} и завершает доказательство  теоремы \ref{th1}.

\end{fulltext}

\end{document}